\documentclass[12pt]{article}
\usepackage{amsmath,amsthm,amscd,amssymb}
\usepackage{latexsym}
\newtheorem{theorem}{Theorem}[section]
\newtheorem{lemma}{Lemma}[section]
\newtheorem{proposition}{Proposition}[section]
\newcommand{\bpropa}{\noindent{\bf Proposition A.}\quad }
\newcommand{\bdefinition}{\noindent{\bf Definition.}\quad }
\newcommand{\bpropb}{\noindent{\bf Proposition B.}\quad }
\newcommand{\bpfprop}{\noindent{\bf Proof of Proposition
2.2.}\quad }
\newcommand{\bpropc}{\noindent{\bf Proposition C.}\quad }

\newtheorem{corollary}{Corollary}[section]

\newcommand{\bproof}{\noindent{\bf Proof: }}
\newcommand{\eproof}{\hfill $\Box$\\}
\newcommand{\bremark}{\noindent{\bf Remark: }}
\newcommand{\eremark}{\hfill \\}

\newcommand{\cplus}{\subset\hspace{-.398cm}\raisebox{.043cm}{$\scriptstyle{+}$}}
\newcommand{\id}{{\rm id}}
\renewcommand{\Re}{{\rm Re}}

\title{Kernels of surjections from ${\cal L}_1$-spaces with
an application to Sidon sets}
\author{N.J.\ Kalton\thanks{Supported in part by NSF grant
DMS-9500125.} \and A.\ Pe\l czy\'nski\thanks{Supported in part
by KBN grant 2P30100406}\date{}}

\begin{document}
\baselineskip20pt
\hoffset-.75in
\voffset-.75in

\maketitle

\begin{abstract} If $Q$ is a surjection from $L^1(\mu)$,
$\mu$ $\sigma$-finite, onto a Banach space containing $c_0$
then (*) $\ker Q$ is uncomplemented in its second dual.  If
$Q$ is a surjection from an ${\cal L}_1$-space onto a Banach
space containing uniformly $\ell_n^\infty$ ($n=1,2,\dots$)
then (**) there exists a bounded linear operator from $\ker
Q$ into a Hilbert space which is not 2-absolutely summing.
Let $S$ be an infinite Sidon set in the dual group $\Gamma$
of a compact abelian group $G$.  Then
$L^1_{\tilde{S}}(G)=\{f\in L^1(G):  \hat{f}(\gamma)=0$ for
$\gamma\in S\}$ satisfies (*) and (**) hence
$L^1_{\tilde{S}}(G)$ is not an ${\cal L}_1$-space and is not
isomorphic to a Banach lattice.  \end{abstract}

\vspace{.3cm}

\noindent{\bf 1991 Mathematical Subject Classification:}
46B03, 43A46.

\vspace{.3cm}

\noindent{\bf Keywords:} Surjection, twisted sum, absolutely
summing operator, Sidon set, Radon Nikodym Property,
Dunford-Pettis Property.

\setcounter{section}{-1}

\section{Introduction} \label{sec-0}

In 1981, J.\ Bourgain \cite{B} solved a long outstanding
problem in Banach space theory by showing the existence of
an uncomplemented subspace of $L^1$ which is isomorphic to
$\ell^1$.  In that paper, he raises the question of whether
it is possible to find an uncomplemented translation
invariant subspace of $L^1(G)$, where $G$ is a compact
abelian group, which is isomorphic to $L^1$.  As a special
case he mentions the question of whether the closed linear
subspace of $L^1$ spanned by the complement of the
Rademachers in the Walsh functions is (a) isomorphic to
$L^1$, or (b) an ${\cal L}_1$-space.  Bourgain \cite[Problem
6]{B} attributes the question to Pisier.  As far as we could
trace, it was also previously considered by Kisliakov and
Zippin.

Suppose $G$ is a compact abelian group and $\Gamma$ its
character group.  For any subset $A$ of $\Gamma$, we define
$L^1_A(G)$ as the closure in $L^1(G)$ of the linear span of
$\{\gamma:  \gamma\in A\}$.  We put
$\widetilde{A}=\Gamma\backslash A$.  Answering in negative
questions (a) and (b), we shall show that if
$S\subset\Gamma$ is an infinite Sidon set, then the
canonical image of $L^1_{\tilde{S}}(G)$ is uncomplemented in
its second dual and it is not an ${\cal L}_1$-space
(Corollary \ref{cor-5.1}).

We approach the problem from a purely Banach space point of
view.  Note that if $S$ is a Sidon set, then the map $Q:
L^1(G)\to c_0(S)$ defined by
$Qf=\{\hat{f}(\gamma)\}_{\gamma\in S}$, where $\hat{f}$
denotes the Fourier transform, is a surjection so that $\ker
Q=L^1_{\tilde{S}}(G)$ is the kernel of a quotient map from
$L^1(G)$ onto a space isomorphic to $c_0(S)$.

We first show (Proposition \ref{prop-2.2}) that if $\mu$ is
a finite measure and $E$ is a Banach space containing an
isomorphic copy of $c_0$, then the canonical image of the
kernel of any surjection $Q$ from $L^1(\mu)$ onto $E$ is
uncomplemented in its second dual; consequently, $\ker Q$ is
non-isomorphic to a Banach lattice.  Our argument depends on
an old lifting principle of Lindenstrauss \cite{L}.

We then turn to part (b) of the question.  Here the idea is
to study subspaces $X$ of an ${\cal L}_1$-space $F$ which
are $GT$-spaces.  (Recall, cf.\ \cite{Pi2}, that a Banach
space $X$ is a $GT$-space if every bounded linear operator
from $X$ into $\ell^2$ is absolutely summing.)  Let $E=F/X$
so that $X$ is the kernel of a quotient map onto $E$.  Then
we show (Theorem \ref{thm-3.1}) that $X$ is a $GT$-space, if
and only if every short exact sequence \[ 0\rightarrow
\ell^2\rightarrow Z \rightarrow E\rightarrow 0 \] splits,
i.e.\ in the language of \cite{KPc} every twisted sum of
$\ell^2$ and $E$ is naturally isomorphic to the Cartesian
product $\ell^2\oplus E$.  This leads us to the general
question of characterizing such Banach spaces $E$.  We show
that if every twisted sum of $\ell^2$ and $E$ splits, then
$E$ (i) fails to have any type $p>1$ (Corollary
\ref{cor-4.1}, cf.\ also \cite{Do2}), and (ii) has cotype
$q<\infty$ (Corollary \ref{cor-4.2}).  In particular, if $E$
contains a subspace isomorphic to $c_0$, then $X$ is not a
$GT$-space and a fortiori, by the Grothendieck Theorem fails
to be an ${\cal L}_1$-space.

Coming back to Sidon sets, we would like to mention that our
techniques do not establish whether, if $S$ is an infinite
Sidon set, the space $L^1_{\tilde{S}}(G)$ can have local
unconditional structure.  (However, see the remark at the
end of the section.)  Furthermore, we do not know whether
the spaces $L^1_{\tilde{S}}(G)$ depend essentially on the
choice of Sidon set (i.e.\ if $S_1$ and $S_2$ are infinite
countable Sidon sets in $\Gamma_1$ and $\Gamma_2$
respectively, are the spaces $L^1_{\tilde S_1}(G_1)$ and
$L^1_{\tilde S_2}(G_2)$ isomorphic?)

In order to keep the paper self contained we include
proofs of several facts on twisted sums which have been known for
about 20 years but which seem to be not available in the
literature.  Many of these facts are contained in the
preprint of Doma\'nski \cite{Do2}.  We are indebted to
Pawe{\l} Doma\'nski, who read the preliminary version of the
paper, for supplying us with additional references and for
many valuable comments.

\noindent{\bf Remark.}\quad After the initial preparation of
the paper, W.B.  Johnson showed, using related techniques,
that the kernel of the quotient map of $L_1$ onto $c_0$
fails to have (Gordon-Lewis) local unconditional structure.

\section{Auxiliary lemmas} \label{sec-1}

In this section we state three essentially known lemmas on
short exact sequences of Banach spaces (cf.\
\cite[Chapt.III]{McL} in the language of homology; \cite{P},
\cite{KPc}, \cite{Ki}, \cite{Do1}, \cite{Vo1}, \cite{V2} in
the setting of Banach spaces and topological vector spaces).

If $u:  X\to Y$ is a (bounded) linear operator acting
between normed spaces $X$ and $Y$, then we put \[
\rho(u)=\inf\{\eta>0:  \forall y\in u(X)\ \exists \ x\in
u^{-1}(y) \ \mbox{with $\|x\|\le\eta\|y\|$}\}.  \]

\begin{lemma} \label{lemma-1.1} Let $E,X,F$ be Banach
spaces, $u:  E\to X$ an isomorphic embedding, $v:  X\to F$ a
surjection, $u(E)=\ker v$.  Then there exists an equivalent
norm on $X$ and $\beta>0$ such that if $X_1$ denotes $X$
equipped with the new norm, then $u:  E\to X_1$ is a linear
isometric embedding and $\beta v:  X_1\to F$ is a quotient
map, hence $E$ is isometrically isomorphic to a subspace of
$X_1$ and $F$ is isometrically isomorphic to the quotient of
$X_1$ by this subspace.  \end{lemma}

\bproof Let $B_Y$ denote the unit ball of a normed space $Y$
and let $cv(W)$ denote the absolute --- for real spaces
(resp.\ circled --- for complex spaces) closed convex hull
of a set $W\subset Y$.

We define the new norm on $X$ to be the gauge functional of
the set \[ cv(u(B_E)\cup \alpha B_X)\cap\beta v^{-1}(B_F) \]
where the positive numbers $\alpha$ and $\beta$ are chosen
so that \[ cv(u(B_E)\cup\alpha B_X) \cap u(E)=u(B_E) \
\mbox{and} \ v(\alpha B_X)\supset\beta B_F.  \] The
existence of $\alpha>0$ and $\beta>0$ in question follows
from the assumptions that $u$ is an isomorphic embedding and
that $v$, being a surjection, is open.  \eproof

Our next lemma in the setting of Banach spaces is often called
``Kisliakov's Lemma''.

\begin{lemma} \label{lemma-1.2} Let $X,Y,X_1$ be Banach
spaces, let $X$ be a subspace of $X_1$ and let $u:  X\to Y$
be a bounded linear operator.  Then there exist a Banach
space $Y_1$ and a linear operator $u_1:  X_1\to Y_1$ such
that $Y$ is a subspace of $Y_1$, $u_1$ is a norm preserving
extension of $u$ and the quotient spaces $X_1/X$ and $Y_1/Y$
are isometrically isomorphic.  Precisely the following
diagram commutes $$\begin{CD}X @>j>> X_1 @>q>> X_1/X\\ @VVuV
@VV{u_1}V @VVIV\\ Y @>>J> Y_1 @>>Q> Y_1/Y\end{CD}$$

\noindent where $j$ and $J$ are natural inclusions, $q$ and
$Q$ quotient maps, and $I$ is an isometric isomorphism.
Moreover, if $u$ is an isometric isomorphism or an isometric
embedding, then so is $u_1$; in general $\rho (u)=\rho
(u_1)$.  \end{lemma}

For a proof except the ``moreover part'' see \cite[pp.\
316--317]{DJTg}.

\vspace{.3cm}

\noindent{\bf Proof of the ``moreover part'':} Without loss
of generality assume that $\|u\|=1$.  Our assumption says
that there is $c$, with $0<c\le 1$ such that $\|u(x)\|_Y\ge
c\|x\|_X$ for $x\in X$.  Recall that $Y_1$ is defined to be
the quotient space of the $\ell^1$-sum $X_1\oplus_1 Y$ by
the subspace \[ W = \{(x,-u(x))\in X_1\oplus_1 Y:  x\in X\},
\] and $u_1$ is defined to be the restriction of the
quotient map $X_1\oplus_1 Y\to Y_1$ to the subspace
$X_1\oplus_1 \{0\}$ naturally identified with $X_1$.  For
fixed $x_1\in X_1$ we have \begin{eqnarray*}
\|u_1(x_1)\|_{Y_1} &=& \inf_{x\in X}
(\|x_1-x\|_{X_1}+\|u(x)\|_Y)\\ & \ge & \inf_{x\in
X}(\|c(x_1-x)\|_{X_1} +c\|x\|_X)\\ & \ge & c\|x_1\|_{X_1}.
\end{eqnarray*} \eproof

We also need a dual version of Kisliakov's Lemma.

\begin{lemma} \label{lemma-1.3} Let $Y,Z,Z_1$ be Banach
spaces.  Let $Z$ be a quotient of $Y$ via the quotient map
$q$ and let $u:  Z_1\to Z$ be a bounded linear operator.
Then there exist a Banach space $Y_1$ such that $Z_1$ is a
quotient of $Y_1$ via the quotient map $q_1$ and a linear
operator $u_1:  Y\to Y_1$ with $\|u\|=\|u_1\|$ such that
$qu_1=uq_1$ and the spaces $X= \ker q$ and $X_1=\ker q_1$
are isometrically isomorphic.  Precisely the following
diagram commutes $$\begin{CD} X @>j>> Y @>q>> Z\\ @AAIA
@AAu_1A @AAuA\\ X_1 @>>J_1> Y_1 @>>q_1> Z_1\end{CD}$$ where
$I$ is an isometric isomorphism and $j$ and $j_1$ are
natural inclusions.

Moreover, if $u$ is a quotient map onto a subspace of $Z$,
then so is $u_1$, precisely $\rho(u)=\rho(u_1)$.
\end{lemma}

\bproof Put \[ Y_1=\{(y,z_1) \in Y\oplus_\infty Z_1:
q(y)=u(z_1)\} \] where the norm in $Y\oplus_\infty Z_1$ is defined
by $\|(y,z_1)\|=\max(\|y\|,\|z_1\|)$.

Define $q_1:  Y_1\to Z_1$ and $u_1:  Y_1\to Y$ by \[
q((y,z_1))=z_1\quad\mbox{and} \ u_1((y,z_1))=y.  \] We omit
the routine verification.  \eproof

\section{Quotient maps from an $L^1$-space whose kernels are
uncomplemented in their second duals} \label{sec-2}

\setcounter{equation}{0}

We begin with a result which was known to several experts in
the field.  It generalizes an old theorem of Lindenstrauss
\cite{L}.  To make the paper self-contained we include the
proof which is essentially the same as Lindenstrauss'
original argument.  For related references cf.\
\cite{KPc} and \cite[Chapt.\ VI]{KpCro} where the results
refer to $p$-homogeneous spaces ($0<p<1$); \cite{Do2} and
\cite{Do3} in the setting of operator ideals.

Here and in the sequel we identify a Banach space with its
canonical image in its second dual.

\begin{proposition} \label{prop-2.1} (Lindenstrauss Lifting
Principle) Let $Y$ and $E$ be Banach spaces and let $Q:
Y\to E$ be a surjection.  Assume \begin{equation}
\label{eq-2.1} \ker Q \quad\mbox{is complemented in its
second dual.} \end{equation}

Then for every ${\cal L}_1$-space $F$ every bounded linear
operator $T:  F\to E$ admits a lifting, i.e.\ there exists a
bounded linear operator $\widetilde{T}:  F\to Y$ such that
$Q\widetilde{T}=T$ \end{proposition}

\bproof The assumption that $F$ is an ${\cal L}_1$-space
means that there exist a $\sigma\in[1,\infty)$ and a subnet
$(F_\alpha)_{\alpha\in\Omega}$ of the net of all finite
dimensional subspaces of $F$ directed by inclusion such that
each $F_\alpha$ is at most $\sigma$ isomorphic to
$\ell^1_{\dim F_\alpha}$.  Let $T_\alpha=T_{|F_\alpha}$,
where $T_{|F_\alpha}$ denotes the restriction of $T$ to
$F_\alpha$.  The lifting property of $\ell^1_{\dim
F_\alpha}$ yields the existence of a linear operator
$\widetilde{T}_\alpha:  F_\alpha\to Y$ with
$Q\widetilde{T}_\alpha = T_\alpha$ and
$\|\widetilde{T}_\alpha\|\le\sigma$.  Since $Q$ is a
surjection, the open mapping theorem yields $\rho(Q)>0$.
Thus, given $\eta>\rho(Q)$, there exists a function
$\varphi:  E\to Y$ (in general neither linear nor
continuous) such that $Q\varphi(e)=e$ and
$\|\varphi(e)\|\le\eta\|e\|$ for $e\in E$.  Now, put
$Y_E=\ker Q$ and, for $0\le r<\infty$ let \[ B(r) =
\{y^{**}\in (Y_E)^{**}:  \|y^{**}\|\le
r\}\quad\mbox{equipped with the $(Y_E)^*$ topology of
$(Y_E)^{**}$.} \] Then $B(r)$ is a compact topological
space.  Hence, by the Tychonoff theorem, the product \[
\prod = \prod_{f\in F} B((\sigma+\eta)\|T\| \|f\|) \] is
also compact.  For every $\alpha\in\Omega$ define
$\pi_\alpha\in\Pi$ by \[ \pi_\alpha(f) =
\left\{\begin{array}{cl} \widetilde{T}_\alpha(f)-\varphi
T(f) & \mbox{for $f\in F_\alpha$} \\0 & \mbox{for $f\notin
F_\alpha$} \end{array}\right.  \] Let $\pi$ be a limit point
of the set $(\pi_\alpha)_{\alpha\in\Omega}$; the existence
of $\pi$ is a consequence of the compactness of $\Pi$.

Let us put \[ \widetilde{T} = P\pi+\varphi T, \] where $P:
(Y_E)^{**}\to Y_E$ is the projection granted by
(\ref{eq-2.1}).

To verify that $\widetilde{T}$ is the desired operator,
first note that $\widetilde{T}_\alpha(f)-\varphi T(f)\in
Y_E$ for every $f\in F_\alpha$ and for $\alpha\in\Omega$;
hence $\pi_\alpha(f)\in B((\sigma+\eta)\|T\| \|f\|)\cap Y_E$
for every $f\in F$ and for every $\alpha\in\Omega$.  Thus
$\pi(f)\in B((\sigma+\eta)\|T\| \|f\|)$ and $P\pi(f)\in
B((\sigma+\eta)\|P\| \|T\| \|f\|)\cap Y_E$.  Since $\|P\|\ge
1$, we get \begin{equation} \label{eq-2.2}
\|\widetilde{T}(f)\|\le \|P\| (\sigma+2\eta)\|T\| \|f\|
\quad\mbox{for $f\in F$.} \end{equation}

To complete the proof one has to verify the linearity of
$\widetilde{T}$.  In view of (\ref{eq-2.1}) and the fact
that $\bigcup_{\alpha\in\Omega}F_\alpha =F$, it is enough to
show \begin{equation} \label{eq-2.3} \widetilde{T}(f'+f'') =
\widetilde{T}(f')+\widetilde{T}(f'') \quad\mbox{for
$f',f''\in\bigcup_{\alpha\in\Omega} F_\alpha$.}
\end{equation} To this end pick $\alpha_0$ so that
$f',f''\in F_{\alpha_0}$ and let $F_\alpha\supset
F_{\alpha_0}$.  Taking into account that
$\widetilde{T}_\alpha(f'+f'')-\widetilde{T}_\alpha(f')-
\widetilde{T}_\alpha(f'') =0$ and $\varphi T(f'+f'')-\varphi
T(f')-\varphi T(f'')\in Y_E$ we get \begin{eqnarray*}
\pi_\alpha(f'+f'')-\pi_\alpha(f')-\pi_\alpha(f'')
&=&-[\varphi T(f'+f'')-\varphi T(f')-\varphi T(f'')]\\ &&
\mbox{for $\alpha\in\Omega$ such that $F_\alpha\supset
F_{\alpha_0}$.} \end{eqnarray*} Thus, remembering that $\pi$
is a limiting point of the net
$(\pi_\alpha)_{\alpha\in\Omega}$, we get \[
\pi(f'+f'')-\pi(f')-\pi(f'') =-[\varphi T(f'+f'')-\varphi
T(f')+\varphi T(f'')]\in Y_E.  \] Applying to both sides of
the latter identity $P$, we get \[
P\pi(f'+f'')-P\pi(f')-P\pi(f'')=-[\varphi T(f'+f'')-\varphi
T(f')-\varphi T(f'')].  \] which yields (\ref{eq-2.3}).
\eproof

Next we discuss relationships of Lindenstrauss' Lifting
Principle with the Radon Nikodym Property (=RNP).  We follow
the terminology and notation of \cite{DU}.

Recall (\cite[Chap.\ III]{DU}) that a linear operator $T:
L^1(\mu)\to E$ ($E$ a Banach space; $\mu$ a finite measure
on a measure space ($\Omega,\Sigma,\mu$)) is representable
if there exists a Bochner integrable function $e(\cdot)\in
L^\infty(\mu;E)$ such that $\|e(\cdot)\|_\infty =\|T\|$ and
\[ Tf=\int_\Omega f(s)e(s)\mu(ds)\quad\mbox{for $f\in
L^1(\mu)$.} \]

A Banach space $E$ has RNP provided for every finite measure
$\mu$ (equivalently for some non purely atomic $\mu$) every
bounded linear operator from $L^1(\mu)$ into $E$ is
representable.

It is interesting to compare Proposition \ref{prop-2.1} with
the well known

\noindent{\bf Fact.} {\em The assertion of Proposition
\ref{prop-2.1} remains valid if the assumption
(\ref{eq-2.1}) is replaced by \begin{equation}
\label{eq-2.4} \mbox{$E$ has RNP and $F$ is an abstract
$L$-space.} \end{equation} }

\bproof Assume first that $F=L^1(\mu)$ with $\mu$ finite.
Then every $e(\cdot)\in L^\infty(\mu;E)$ can be represented
as a sum of an absolutely convergent series in
$L^\infty(\mu;E)$ of countably valued functions.  Hence, by
the open mapping principle there exist $\delta>0$ and a
function $y(\cdot)\in L^\infty(\mu;Y)$ such that
$Q(y(s))=e(s)$ for $s\in \Omega$ $\mu$ a.e.\ and
$\|y(\cdot)\|_\infty\le\delta\|e(\cdot)\|_\infty$.  Thus if
$e(\cdot)$ represents $T:  L^1(\mu)\to E$, then we define
$\widetilde{T}:  L^1(\mu)\to Y$ by \[ \widetilde{T}f=\int
f(s)y(s)\mu(ds)\quad\mbox{for $f\in L^1(\mu)$.} \]

The general case follows from the observation that by
Kakutani's representation theorem every abstract $L$-space
is the $\ell^1$-sum of a family of spaces
$(L^1(\mu_\alpha))_{\alpha\in A}$ with $\mu_\alpha$ finite
for all $\alpha\in A$.  \eproof

A simple consequence of Proposition \ref{prop-2.1} is

\begin{corollary} \label{cor-2.1} If $E$ fails RNP and $Q:
\ell^1(A)\to E$ is a surjection, then $\ker Q$ is not
complemented in $(\ker Q)^{**}$.  \end{corollary}

\bproof If $E$ fails RNP, then there exists a finite measure
$\mu$ and a bounded linear operator $T:  L^1(\mu)\to E$
which is not representable.  Since every bounded linear
operator $u:  L^1(\mu)\to \ell^1(A)$ is representable
(\cite[p.83, Corollary 8]{DU}), so is $Qu$.  If $\ker Q$
were complemented in $(\ker Q)^{**}$ then, by Proposition
\ref{prop-2.1}, $T$ would be of the form $Qu$, a
contradiction.  \eproof

Our next result is more sophisticated

\begin{proposition} \label{prop-2.2}

Let $Q:  L^1(\mu)\to E$ be a surjection.  Assume that the
measure $\mu$ is finite and $E$ contains a subspace, say
$E_1$, isomorphic to $c_0$.  Then $\ker Q$ is uncomplemented
in $(\ker Q)^{**}$.  \end{proposition}

In particular the kernel of a surjection of $L^1(\mu)$ onto
$c_0(A)$ is uncomplemented in its second dual whenever $A$
is infinite and $\mu$ is a finite measure.

\bproof Assume first that $E$ is separable.  Then, by
Sobczyk's Theorem (\cite[I, Theorem 2.f.5]{LT}), there
exists a projection $P:  E
\stackrel{\mbox{\scriptsize{onto}}}{\longrightarrow} E_1$.
Let $(e_n,e^*_n)^\infty_{n=1}$ be the biorthogonal system in
$(E_1,E_1^*)$ induced by the unit vector basis of $c_0$.
Put $\varphi_n=(PQ)^*(e_n^*)$ for $n=1,2,\dots$.  Then
$(\varphi_n)\subset L^\infty(\mu)=[L^1(\mu)]^*$.  Regarding
$(\varphi_n)$ as a sequence in $L^2(\mu)$ we infer that
$\varphi_n\to 0$ weakly in $L^2(\mu)$ as $n\to\infty$
(because $e_n^*\to 0$ in the $c_0$ topology of
$\ell^1=(c_0)^*$ as $n\to 0$).  By Mazur's Theorem some
convex combinations of the $\varphi_n$'s tend to 0 strongly
in $L^2(\mu)$.  Hence, by a result of F.\ Riesz, a
subsequence of these convex combinations tends to 0
$\mu$-almost everywhere.  Thus there is an increasing
sequence of the indices $0=k_0<k_1<\cdots$, a sequence
$(\psi_n)\subset L^\infty(\mu)$ such that
$\psi_n=\sum^{k_n}_{j=k_{n-1}+1} a_j\varphi_j$ with $a_j\ge 0$
and $\sum^{k_n}_{j=k_{n-1}+1} a_j=1$ ($n=1,2,\dots$), and
\begin{equation} \label{eq-2.5} \lim_n
\psi_n(s)=0\quad\mbox{for $s\in\Omega$ $\mu$-a.e.}
\end{equation} We put \[ R\left(\sum^\infty_{j=1}
t_je_j\right) = \sum_{n=1}^\infty
\left(\sum^{k_n}_{j=k_{n-1}+1}
t_ja_j\right)\left(\sum_{j=k_{n-1}+1}^{k_n} e_j\right)
\mbox{for $(t_j)\in c_0$.} \] Then $R$ is a projection from
$E_1$ onto its subspace $E_0$ isomorphic to $c_0$ and
spanned by the sequence of ``characteristic functions'',
$\left(\sum^{k_n}_{j=k_{n-1}+1} e_j\right)^\infty_{n=1}$.
Thus the natural embedding $J:  E_0\to E$ satisfies $RPJ=
\id_{E_0}$.  Clearly $E_0$ being isomorphic to $c_0$ fails
RNP.  Thus there is a bounded linear operator $T:  L^1\to
E_0$ which is not representable ($L^1$ denotes the space of
absolutely Lebesgue integrable scalar valued functions on
$[0,1]$).  Now, if $\ker Q$ were complemented in $(\ker
Q)^{**}$ then, by Proposition \ref{prop-2.1}, there would
exist a bounded linear operator $\widetilde{S}:  L^1\to
L^1(\mu)$ such that $Q\widetilde{S}=JT$.  Thus
$RPQ\widetilde{S}=RPJT=T$.  Now, observe that \[ RPQ(f) =
\sum^\infty_{n=1} \int_\Omega
f(s)\psi_n(s)\mu(ds)\left(\sum_{j=k_{n-1}+1}^{k_n}
e_j\right)\quad\mbox{for $f\in L^1(\mu)$.} \] Note that the
sequence $(\sum^{k_n}_{j=k_{n-1}+1}e_j)^\infty_{n=1}$ is
equivalent to the unit vector basis of $c_0$.  Thus the
condition (\ref{eq-2.5}) yields that the operator RPQ is
representable (cf.\ \cite[p.75, Remark after the proof of
Lemma 2.11]{DU}).  Hence $T=RPQ\widetilde{S}$ would be
representable because an operator from $L^1(\mu)$ to a
Banach space is representable iff it factors through
$\ell^1(A)$ (cf.\ \cite[Chapt.\ III, \S1, proof of Theorem
8]{DU}), a contradiction.  \eproof

The argument for non separable $E$ is almost the same.
Instead of Sobczyk's Theorem we use the following
generalization.

\begin{lemma} \label{lemma-2.1} Let $E$ be a quotient of
$L^1(\mu)$ with $\mu$-finite measure.  Assume that $E$
contains a subspace $E_1$ isomorphic to $c_0$.  Then $E_1$
is complemented in $E$.  \end{lemma}

\bproof If $\mu$ is finite then the natural injection of
$L^2(\mu)$ into $L^1(\mu)$ is bounded and has a dense range,
hence $L^1(\mu)$ is a WCG space.  Thus $E$ is a WCG space.
Therefore every separable subspace of $E$ is contained in a
separable subspace which is a range of a contractive
projection from $E$ (cf.\ \cite[pp.\ 237--240]{DeGoZ}).
Combining this fact with Sobczyk's Theorem we get the
desired conclusion.  \eproof\newpage

\noindent{\bf Remarks} \begin{enumerate} \item The assertion
of Proposition \ref{prop-2.2} remains valid if the
assumption ``$\mu$-finite'' is replaced by ``$\mu$ $\sigma$-
finite'' because every $L^1(\nu)$ with $\nu$ $\sigma$-finite is
isomorphic as a Banach space with $L^1(\mu)$ for some finite
measure $\mu$.  \item After reading a preliminary version of
this paper, S.\ Kwapien has shown us an alternative proof of
Proposition \ref{prop-2.2} which does not use the Lindenstrauss
Lifting Principle.
We present his argument with his permission.
\end{enumerate}

Let $X_0$ be a subspace of a Banach space $X$.  Then
$(X/X_0)^*$ can be identified with $X_0^\bot$ and $X^{**}_0$
with $X_0^{\bot \bot}$ where \begin{align*}X^\bot_0 & =
\{x^*\in X^*:x^*(x)=0\text{ for } x\in X_0\},\\ X^{\bot
\bot}_0 & = \{x^{**}\in X^{**}:x^{**}(x^{*})=0\text{ for
}x^*\in X^\bot_0\}.\end{align*} The subspace $X+X^{\bot
\bot}_0$ of $X^{**}$ is norm closed.  The condition $X_0$ is
complemented in $X^{**}_0$ is equivalent to the existence of
a bounded linear projection $p:X_0^{\bot \bot
}\underset{\text{onto}}{\longrightarrow}X_0$.

For $j=0,1,2$ denote by $(\delta ^{(j)}_n)_{n=1}^\infty $
the unit vector basis of $c_0$, $\ell^1$, $\ell^2$
respectively.

\bpropa {\em If $X_0$ is complemented in $X_0^{**}$ then for
every bounded operator $T:c_0\to X/X_0$ there exists a
weakly null sequence $(x_n)\subset X$ such that
$T(\delta^{(0)}_n)=Q(x_n)$ for $n=1,2,\ldots $ where $Q:X\to
X/X_0$ is the quotient map.}

\bproof First note that the formula
$$\tilde{p}(x+x^{**})=x+p(x^{**})\quad (x\in X,x^{**}\in
X_0^{\bot \bot})$$ well defines a projection from $X\oplus X_0^{\bot
\bot}$ onto
$X$ with
$\| \tilde{p}\|\leq \|p\|+2$.

Let $S=IT^*:X^\bot_o\overset{T^*}{\longrightarrow}\ell^1
\overset{I}{\longrightarrow}\ell^2$ where $I:\ell^1\to
\ell^2$ is the natural embedding since $I$ is $2$-summing,
so is $S$.  Thus $S$ extends to a bounded operator
$\tilde{S}:X^*\to \ell^2$.  Put
$x^{**}_n=(\tilde{S})^*(\delta^{(2)}_n)$ for $n=1,2,\ldots
$. Then $x^{**}_n\in X\in X_0^{\bot \bot}$.  Indeed, pick
$y_n\in X$ so that $Q(y_n)=T(\delta_n^{(0)})$.  Then
$x_n^{**}=y_n+(x^{**}_n-y_n)$ and $x^{**}_n-y_n\in X_0^{\bot
\bot }$ because for every $x^*\in X^\bot_0$ one has
$$x^{**}_n(x^*)=\delta _n^{(2)}(\tilde{S}x^*)=\delta
^{(2)}_n(Sx^*)=\delta ^{(1)}_n(T^*x^*)=x^*(T(\delta
_n^{(0)})).$$

Now put $x_n=\tilde{p}(x^{**}_n)$ for $n=1,2,\ldots $. Since
$x_n=\tilde{p}(\tilde{S})^*(\delta ^{(2)}_n)$ and $(\delta
^{(2)}_n)$ is a weakly null sequence in $\ell^2$, so is
$(x_n)$ in $X$.  Finally $Q(x_n)=Q(y_n)$ because
$x_n=\tilde{p}(y_n+(x_n^*-y_n))=y_n+p(x^*_n-y_n)$ and
$p(x^{**}_n-y_n)\in X_0$.\eproof

\bdefinition A Banach space $X$ has property $(K)$ if for an
arbitrary weak* null sequence $(\varphi^*_n)\subset X^*$
there exists a CCC sequence $(\psi^*_n)$ such that
\begin{equation}\label{eq2.6} \lim_{k}\psi^*_k(x_k)=0
\text{ for every weakly null sequence }(x_k)\subset
X.\end{equation} ``CCC'' stands for ``consecutive convex
combinations''; $(\psi_k^*)$ is a CCC sequence for
$(\varphi^*_n)$ if there exist an increasing sequence of the
indices $0=n_0< n_1< \ldots $ and a sequence $(\lambda_n)$
of non-negative scalars such that
$$\psi_k^*=\sum^{n_k}_{j=n_{k-1}+1}\lambda_j\varphi_j^*
\text{ with } \sum^{n_k}_{j=n_{k-1}+1}\lambda_j=1\quad
(k=1,2,\ldots ).$$

\bpropb {\em $L^1(\mu )$ has $(K)$ for every finite measure
$\mu $.}

\bproof If $(\varphi^*_n)\subset L^\infty(\mu
)=(L^1(\mu))^*$ is a weak* star null sequence then some
subsequence of $(\varphi^*_n)$ is a weakly null sequence in
$L^2(\mu )$.  Hence, by Mazur's theorem, some CCC sequence
$(\psi_k^*)$ for $(\varphi^*_n)$ strongly converges to zero
in $L^2(\mu )$.  Since weakly null sequences in $L^1$ are
equi-integrable one easily
checks that $(\psi^*_k)$ satisfies \eqref{eq2.6} (cf. the
argument at the end of the proof of Lemma 5.2)\eproof

The space $c_0$ fails $(K)$ in the following strong sense:

\bpropc {\em If $(y_k^*)$ is a CCC sequence for
$(\delta^{(1)}_n)\subset \ell^1=(c_0)^*$ then there exists a
bounded linear projection $V:c_0\to c_0$ such that
\begin{equation}\label{eq2.7} \inf_k
y^*_k(V(\delta^{(0)}_k))>0\quad \Box \end{equation}}

\bpfprop Let $U:c_0\to E$ be an isomorphic embedding.  Since
Sobczyk's theorem generalizes to quotients of $L^1(\mu )$
for $\mu $ finite (Lemma 2.1), there exists a projection,
$P:E\to E$ with $P(E)=U(c_0)$.  Put
$\varphi^*_n=Q^*P^*(U^{-1})^*(\delta ^{(1)}_n)\in (L^1(\mu
))^*$ for $n=1,2,\ldots$.  Since $(\delta ^{(1)}_n)$ is a
weak* null sequence in $\ell^1=(c_0)^*$, so is
$(\varphi^*_n)$ in $(L^1(\mu ))^*$.  By Proposition B there
exists a CCC sequence $(\psi^*_k)$ for $(\varphi ^*_n)$
which satisfies \eqref{eq2.6}.  Let $(y^*_k)$ be the CCC
sequence for $(\delta ^{(1)}_n)$ of the same convex
combinations as the sequence $(\psi^*_k)$ for
$(\varphi^*_n)$.  Let $V:c_0\to c_0$ be that of Proposition
C. Put $T:UV:c_0\to E$.  Assume that $\ker Q=X_0\subset
L^1(\mu )$ were complemented in $X^{**}_0$.  Then, by
Proposition A (applied to $X=L^1(\mu ))$ there would exist a
weakly null sequence sequence $(x_k)\subset L^1(\mu )$ such
that $T(\delta ^{(0)}_k)=Q(x_k)$ for $k=1,2,\ldots$.  On the
other hand we have $$\psi ^*_k(x_k)=y^*_k(V(\delta
^{(0)}_k)\quad \text{for }k=1,2,\ldots $$ This would lead to
a contradiction with \eqref{eq2.7}, because
$\lim_k\psi_k^*(x_k)=0$ by \eqref{eq2.6}.\eproof

Finally note that Proposition B and accordingly Proposition
2.2 can be generalized to preduals of $C^*$-algebras with
finite faithful trace.

Our last result in this section is a partial converse to
Corollary \ref{cor-2.1}.

\begin{proposition} \label{prop-2.3} Suppose $E$ has RNP and
is complemented in $E^{**}$.  Then for every quotient map
$Q:  \ell^1(A)\to E$, $\ker Q$ is complemented in $(\ker
Q)^{**}$.  \end{proposition}

\bproof Denote by $P$ a bounded projection from $E^{**}$
onto $E$ and by $\Pi:  \ell^1(A)^{**}\to\ell^1(A)$ the
canonical projection.

We start by observing that if $Z$ is an abstract $L$-space,
then every operator $T:  Z\to E$ factors through a space
$\ell^1(B)$.  This follows from the fact that $Z$ can be
decomposed as an $\ell^1$-sum of $L^1(\mu_\alpha)$ spaces
where each $\mu_\alpha$ is a finite measure and from the
Lewis-Stegall theorem (cf.\ \cite{LT}, \cite[Chapt.\ III,
\S1, Theorem 8]{DU}).  Since $\ell^1(B)$ is projective, this
implies that $T$ also factors through $Q$.  We apply these
remarks to $Z=(\ell^1(A))^{**}$ and to $T=PQ^{**}$ to deduce
the existence of a bounded linear operator $S:
[\ell^1(A)]^{**}\to \ell^1(A)$ such that $PQ^{**}=QS$.

Now let $V=S+(I-S)\Pi:  [\ell^1(A)]^{**}\to\ell^1(A)$.
Clearly $V$ is a projection onto $\ell^1(A)$.  Furthermore
$V(\ker Q^{**})=\ker Q$ (indeed let $x^{**}\in\ker Q^{**}$.
Put $x=\pi x^{**}\in\ell^1(A)$.  Then
$QSx^{**}=PQ^{**}x^{**}=0$ and $PQ^{**}x=PQx=Qx$.  Therefore
$QVx^{**}=Q(I-S)x=$ $Qx-PQ^{**}x=Q$).  However $\ker Q^{**}$
is the weak*-closure of $\ker Q$ in $[\ell^1(A)]^{**}$
(since $Q$ is a quotient map) and so is naturally isomorphic
to $(\ker Q)^{**}$.  \eproof

\section{Subspaces of an ${\cal L}_1$-space which are GT
spaces} \label{sec-3}

\setcounter{equation}{0}

Recall that a Banach space $Y$ is a twisted sum of Banach
spaces $X$ and $Z$, in symbols $Y=X \cplus \ Z$ provided \[
0\rightarrow X \stackrel{j}{\rightarrow} Y
\stackrel{q}{\rightarrow} Z \rightarrow 0 \] is a short
exact sequence, i.e.\ $j(X)=\ker q$, with $j$ being an
isometrically isomorphic embedding and $q$ being a quotient
map.  We say that a twisted sum $Y=X\cplus \ Z$ splits
provided it is naturally isomorphic to the Cartesian product
$X\oplus Z$, i.e.\ there exists a bounded linear operator
$v:  Z\to Y$ such that \begin{equation} \label{eq-3.1} qv =
\id_Z \end{equation} where $\id_Z$ denotes the identity
operator on $Z$.  Note that if $v$ satisfies (\ref{eq-3.1})
then $\id_Y-vq$ is a projection from $Y$ onto $X$.
Conversely, if $p:  Y\to X$ is a bounded linear projection
onto $X$ then the formula $v(z)=y-p(y)$ for any $y\in Y$
such that $q(y)=z$ well defines $v:  Z\to Y$ which satisfies
(\ref{eq-3.1}), hence $X\cplus \ Z$ splits.

Next recall that every Banach space is a quotient space of
$\ell^1(A)$ for an appropriate set $A$.

We begin with a ``purely formal'' but useful fact.

\begin{proposition} \label{prop-3.1} Suppose we are given
Banach spaces $E$ and $H$ and a quotient map $q_E:
\ell^1(A)\to E$ and let $X_E=\ker q_E$.

Then the following conditions are equivalent \begin{itemize}
\item[(i)] every bounded linear operator from $X_E$ into $H$
extends to a bounded linear operator from $\ell^1(A)$ into
$H$; \item[(ii)] every twisted sum $H\cplus \ E$ splits.
\end{itemize} \end{proposition}

\bproof (ii) $\Rightarrow$ (i).  By Lemma \ref{lemma-1.2} a
bounded linear operator $u:  X_E\to H$ extends to a linear
operator $u_1:  \ell^1(A)\to H\cplus \ E$.  By (ii) there
exists $v:  E\to H\cplus \ E$ with $qv=\id_E$ where $q:
H\cplus \ E\to E$ is the quotient map with $\ker q = H$.
Put $p= \id_{H
{\subset\hspace{-.2cm}\raisebox{.02cm}{$\scriptscriptstyle{+}$}}
E}-vq$.  Then $p$ is a projection from $H\cplus \ E$ onto
$H$ and $pu_1$ is the desired extension of $u$.

(i) $\Rightarrow$ (ii).  Fix a twisted sum $Y=H
{\subset\hspace{-.408cm}\raisebox{.043cm}{$\scriptstyle{+}$}}
\ E$ and let $q:  Y\to E$ be the quotient map with $\ker
q=H$.  Fix $c>1$.  The lifting property of $\ell^1(A)$
yields the existence of a bounded linear operator $\varphi:
\ell^1(A)\to Y$ such that $q\varphi=q_E$ and
$\|\varphi\|<c$.  Obviously $\varphi(X_E)\subset H=\ker q$.
Thus, by (i), the restriction of $\varphi$ to $X_E$ extends
to a bounded linear operator, say $v:  \ell^1(A)\to H$.  Let
us consider the operator $\varphi-v:  L^1(A)\to Y$.  Clearly
$\ker(\varphi -v)\supset X_E$.  Hence $\varphi-v$ factors
through $q_E$, precisely the formula $u(e)=(\varphi-v)(\xi)$
for $e=q_E(\xi)\in E$ well defines a bounded linear operator
$u:  E\to Y$ such that $uq_E=\varphi-v$.  Since $qv=0$, we
have \[ quq_E= q\varphi =q_E.  \] Thus $qu=\id_E$ because
$q_E(\ell^1(A))=E$.  Hence $H\cplus \ E$ splits.  \eproof

\bremark For an analogue of Proposition \ref{prop-3.1} in
Frechet spaces cf.\ \cite{V2}, \cite{V3}.  \eremark

Under the additional assumption that $H$ is complemented in
its second dual, Proposition \ref{prop-3.1} generalizes to
${\cal L}_1$-spaces.

\begin{proposition} \label{prop-3.2} Let $E$ and $H$ be
Banach spaces.  Assume that $H$ is complemented in its
second dual.  Then the following conditions are equivalent
\begin{itemize} \item[(ii)] every twisted sum $H\cplus \ E$
splits; \item[(iii)] for every ${\cal L}_1$-space $F$ and
every quotient map $q_E:  F\to E$ every bounded linear
operator from $X_E=\ker q_E$ into $H$ extends to a bounded
linear operator from $F$ into $H$; \item[(iv)] there exists
an ${\cal L}_1$-space $F$ and a quotient map $q_E:  F\to E$
such that $X_E=\ker q_E$ has the above extension property.
\end{itemize} \end{proposition}

\bproof (ii) $\Rightarrow$ (iii).  The proof is the same as
that of the implication (ii) $\Rightarrow$ (i) of
Proposition \ref{prop-3.1}.

(iii) $\Rightarrow$ (iv).  Trivial.

(iv) $\Rightarrow$ (ii).  The proof differs from the proof
of the implication (i) $\Rightarrow$ (ii) of Proposition
\ref{prop-3.1} only how the bounded linear operator
$\varphi$ which lifts $q_E$ is constructed.  Instead of
using the lifting property of $\ell^1(A)$ we apply
Proposition \ref{prop-2.1}.  At this place the assumption
that $H$ is complemented in $H^{**}$ is used.  \eproof

\bremark Condition (ii) has isomorphic character in the
following sense.  If $(H,E)$ is a pair of Banach spaces such
that every twisted sum $H\cplus \ E$ splits and if
$(H_1,E_1)$ is another pair such that $H$ is isomorphic to
$H_1$ and $E$ is isomorphic to $E_1$ then every twisted sum
$H_1\cplus \ E_1$ splits.  This is an immediate consequence
of Lemma \ref{lemma-1.1}.  Thus Propositions \ref{prop-3.1}
and \ref{prop-3.2} remain valid if one replaces $E$ by a
Banach space isomorphic to $E$ or equivalently if $q_E$ is
an arbitrary surjection onto $E$.  \eremark

Next we discuss a qualitative and a local version of
condition (ii).

For a twisted sum $Y=X\cplus \ Z$ put \begin{equation}
\label{eq-3.2} spl(Y) = \left\{\begin{array}{l} +\infty \
\mbox{if $Y$ does not split} \\
\inf(\|v\|^2+\|\id_Y-vq\|^2)^\frac12 \end{array} \right.
\end{equation} where the infimum extends over all $v:  Z\to
Y$ satisfying (\ref{eq-3.1}).

Our next proposition can be deduced from some results in
Doma\'nski's Ph.D.\ Thesis (Pozna\'n 1986) which are stated
in terms of operator ideals (cf.\ \cite[Theorem II.1.1 and
Theorem II.3.1]{Do2}; cf.\ also \cite{Do3}).

\begin{proposition} \label{prop-3.3} For every pair of
Banach spaces $E$ and $H$ the following conditions are
equivalent \begin{itemize} \item[(ii)] every twisted sum
$H\cplus \ E$ splits; \item[(v)] $\sup spl (H\cplus \
E)=C<+\infty$, where the supremum extends over all twisted
sums $H\cplus \ E$.  \end{itemize}

Moreover, if $H$ is complemented in $H^{**}$ and there
exists a family $(H_\alpha)_{\alpha\in\Omega}$ of finite
dimensional subspaces of $H$ directed by inclusion and such
that $\cup_{\alpha\in\Omega}H_\alpha$ is dense in $H$ and
each $H_\alpha$ is the range of projection $\pi_\alpha:
H\to H_\alpha$ with $\sup_{\alpha\in\Omega}
\|\pi_\alpha\|<+\infty$ then the equivalent conditions (ii)
and (v) are equivalent to \begin{itemize} \item[(vi)] $\sup
spl(H_\alpha\cplus \ E) = C_1<\infty$, where the supremum
extends over all $\alpha\in\Omega$ and all twisted sums
$H_\alpha\cplus \ E$.  \end{itemize} \end{proposition}

\bproof (ii) $\Rightarrow$ (v).  Let $L(X,Y)$ denote the
Banach space of all bounded linear operators from $X$ into
$Y$.  A restatement of condition (i) of Proposition
\ref{prop-3.1} says that the restriction operator maps
$L(\ell^1(A),H)$ onto $L(X_E,H)$.  Hence, by the open
mapping principle, there exists an $M\in (0,\infty)$ such
that every $u\in L(X_E,H)$ extends to an $u_1\in
L(\ell^1(A),H)$ with $\|u_1\|\le M\|u\|$.  Now the analysis
of the proof of the implication (i) $\Rightarrow$ (ii)
yields (v) with $C\le \sqrt{M^2+(M+1)^2}$.

(v) $\Rightarrow$ (ii).  Trivial.

(v) $\Rightarrow$ (vi).  Fix a twisted sum $H_\alpha\cplus \
E$ and denote by $i:  H_\alpha\to H$ the natural
inclu\-sion.  By Lemma \ref{lemma-1.2} there exists a
twisted sum $H\cplus \ E$ and an isometrically isomorphic
embedding $I:  H_\alpha\cplus \ E\to H\cplus \ E$ which
extends $i$.  By (v) for every $a>1$ there exists a
projection $p:  H\cplus \ E
\stackrel{\mbox{onto}}{\rightarrow} H$ with $\|p\|\le aC$.
Put $P=\pi_\alpha pI$.  Then $P:  H_\alpha\in E\to H_\alpha$
is a projection with $P(H_\alpha\cplus \ E)=H_\alpha$ and
$\|P\|\le aC\sup_\alpha \|\pi_\alpha\|=C_2$.  Thus
$spl(H_\alpha\cplus \ E)\le C_1$ where
$C_1=\sqrt{C_2^2+(1+C_2)^2}$.

(vi) $\Rightarrow$ (v).  Let
$H{\subset\hspace{-.418cm}\raisebox{.043cm}{$\scriptstyle{+}$}}
\ E$ be a twisted sum.  By Lemma \ref{lemma-1.2} for each
$\alpha\in\Omega$ there exist a twisted sum
$H_\alpha{\subset\hspace{-.418cm}\raisebox{.043cm}{$\scriptstyle{+}$}}
\ E$ and a linear operator $\Pi_\alpha$ which extends
$\pi_\alpha$ and satisfies $\|\Pi_\alpha\|=\|\pi_\alpha\|$.
Fix $a>1$.  By (vi) there exists a projection $p_\alpha:
H_\alpha{\subset\hspace{-.418cm}\raisebox{.043cm}{$\scriptstyle{+}$}}
\ E\to H_\alpha$ with
$p_\alpha(H_\alpha{\subset\hspace{-.398cm}\raisebox{.043cm}{$\scriptstyle{+}$}}
E)=H_\alpha$ and $\|p_\alpha\|\le a C_1$.  Consider the
family of operators $(i_\alpha
p_\alpha\Pi_\alpha)_{\alpha\in\Omega}$.  Clearly \[ i_\alpha
p_\alpha\Pi_\alpha:  H\cplus \ E\to H_\alpha\subset H\subset
H^{**}\ \mbox{and} \ \|i_\alpha p_\alpha \Pi_\alpha\|\le
aC_1\sup_\alpha\|\pi_\alpha\|.  \] Note that
\begin{equation} \label{eq-3.3} \mbox{if $h\in H_\alpha$ and
$H_{\alpha'}\supset H_\alpha$ then $i_\alpha p_\alpha
\Pi_\alpha(h)=i_{\alpha'}p_{\alpha'}\Pi_{\alpha'}(h)=h$.}
\end{equation}

Now, in the Stone-\v{C}ech compactification $\beta(\Omega)$
of the discrete set $\Omega$, consider the family $(cl {\cal
O}_\alpha)_{\alpha\in\Omega}$ where \[ {\cal O}_\alpha =
\{\alpha'\in \Omega:  H_{\alpha'}\supset
H_\alpha\}\quad\mbox{for $\alpha\in\Omega$} \] and $cl W$
denotes the closure of a set $W$.  Clearly the family
$({\cal O}_\alpha)_{\alpha\in\Omega}$ is centered.  Hence
the intersection $\bigcap_{\alpha\in\Omega} cl {\cal
O}_\alpha$ is non-empty.  Pick a $\phi\in
\bigcap_{\alpha\in\Omega} cl {\cal O}_\alpha$.  Denote by
$\lim_\phi(f)$ the evaluation at $\phi$ of the unique
continuous extension of a bounded scalar-valued function $f$
on $\Omega$ to a continuous function on $\beta(\Omega)$.
For each $y\in H\cplus \ E$ and each $h^*\in H^*$ let
$f_{y,h^*}$ be the scalar valued function on $\Omega$
defined by $f_{y,h^*}(\alpha)=[i_\alpha
p_\alpha\prod_\alpha(y)](h^*)$ (we identify $H$ with its
canonical image in $H^{**}$).  For each $y\in H\cplus \ E$
define the function $T(y)$ on $H^*$ by $T(y)(h^*)=\lim_\phi
f_{y,h^*}$.  It can be easily verified that $T(y)\in H^{**}$
and $T:  H\cplus \ E\to H^{**}$ defined by $y\to T(y)$ for
$y\in H\cplus \ E$ is a linear operator with $\|T\|\le a C_1
\sup_\alpha \|\pi_\alpha\|$.  It follows from (\ref{eq-3.3})
and the density of $\bigcap_{\alpha\in\Omega} H_\alpha$ in
$H$ that $T(h)=h$ for $h\in H$.  Now if $S$ is a bounded
linear projection from $H^{**}$ onto $H$ then ${P}=ST$ is
the desired projection from $H\cplus \ E$ onto $H$ with
$\|P\|\le a C_1\|S\|\sup_\alpha \|\pi_\alpha\|=u C_0$.  This
yields (v) with $C\le \sqrt{C_0^2+(C_0+1)^2}$.  \eproof

We get an important application of Proposition
\ref{prop-3.1}--\ref{prop-3.3} by specifying $H$ to be an
infinite dimensional Hilbert space, say $H=\ell^2$.  The
theory of absolutely summing operators enters.  First we
recall some results on absolutely summing operators
essentially due to Grothendieck \cite{Gr1} with Maurey's
\cite{Ma1} improvement of (jjj) for $p<1$ (cf.\ \cite[pp.\
60]{Pi2}).

\begin{itemize} \item[{\bf (G)}] {\em Let $X$ be a closed
linear subspace of an ${\cal L}_1$-space $F$.  Then for
every bounded linear operator $u:  X\to\ell^2$ the following
conditions are equivalent \begin{itemize} \item[(j)] $u$
extends to a bounded linear operator from $F$ into $\ell^2$;
\item[(jj)] $u$ is 2-absolutely summing; \item[(jjj)] $u$ is
$p$-absolutely summing for $p\in [0;2]$; \item[(jjjj)] for
every bounded linear operator $v:  \ell^2\to X$ the
composition $uv$ is in the Hilbert-Schmidt class.
\end{itemize}} \end{itemize}

Recall that a Banach space $X$ is called a GT-space (cf.\
\cite{LT}) if every bounded linear operator from $X$ into
$\ell^2$ is 1-absolutely summing.

Combining (G) with Propositions
\ref{prop-3.1}--\ref{prop-3.3} and with the Remark after
Proposition \ref{prop-3.2} we get

\begin{theorem} \label{thm-3.1} Let $E$ be a Banach space
and let $Q$ be a linear surjection from an ${\cal
L}_1$-space $F$ onto $E$.  Let $X_E=\ker Q$.  Then the
following conditions are equivalent \begin{itemize}
\item[(+)] $X_E$ is a GT-space; \item[(++)] every twisted
sum $\ell^2\cplus \ E$ splits; \item[(+++)] $\sup_n
spl(\ell_n^2\cplus \ E) < +\infty$, where the supremum
extends over all twisted sums $\ell^2_n\cplus \ E$ and over
positive integers $n$.  \end{itemize} \end{theorem}

\section{Banach spaces $E$ with a non-trivial twisted sum of
$\ell^2$ and $E$} \label{sec-4}

\setcounter{equation}{0}

We begin with two known lemmas (cf.\ e.g.\ \cite[Chapt.I \S5
and Chapt.II \S5]{Do2}).

\begin{lemma} \label{lemma-4.1} Let $E,F,H$ be Banach
spaces.  Assume that
there
exist linear operators $\varphi:  E\to F$ and $\psi:  F\to
E$ such that $\varphi\psi = \id_F,\quad \|\psi\|\le C$ and $\
\|\varphi\|\le C,$ so that $F$ is $C^2-$equivalent to a $C$-complemented
subspace of $E.$ Then there exists a twisted sum $H\cplus \ E$ such that
\begin{equation} \label{eq-4.2} spl(H\cplus \ E) \ge C^{-2}
\left(\frac{1}{8} \sup spl \, Y-\frac14\right)
\end{equation} where the supremum extends over all twisted
sums $Y=H\cplus \ F$.

In particular, if there exists a twisted sum $H\cplus \ F$
which does not split, then there is a twisted sum $H\cplus \
E$ which does not split.  \end{lemma}

\bproof Pick a twisted sum $H\cplus \ F$ so that $2
spl(H\cplus \ F) \ge \sup spl(Y)$.  By Lemma \ref{lemma-1.3}
there exist a twisted sum $H\cplus \ E$ and a linear
operator $\Phi:  H\cplus \ E\to H\cplus \ F$ such that
$\|\Phi\|=\|\varphi\|$ and $q\Phi=\varphi q_1$ where $q:
H\cplus \ F\to F$ and $q_1:  H\cplus \ E\to E$ are quotient
maps with $\ker q = \ker q_1=H$.  If $H\cplus \ E$ does not
split, we are done.  Otherwise there exists a bounded linear
operator $v_1:  E\to H\cplus \ E$ such that $q_1v_1=\id_E$
and $\|v_1\|\le 2 spl(H\cplus \ E)$.  Put $v=\Phi v_1\psi:
F\to H\cplus \ F$.  Clearly, \[ qv=q\Phi v_1\psi = \varphi
q_1v_1\psi = \varphi\psi = \id_F \] and \[ \|v\|\le
C^2\|v_1\|\le 2C^2 spl(H\cplus \ E).  \] On the other hand,
\[ 2\|v\|+1\ge spl(H\cplus\; F)\ge 2^{-1}\sup spl(Y).  \]
The last two chains of inequalities obviously yield
(\ref{eq-4.2}).  \eproof

Remembering that the adjoint of a linear isometric embedding
is a quotient map and vice versa, and applying directly
formula (\ref{eq-3.2}) we get

\begin{lemma} \label{lemma-4.2} Let $H\cplus \ E$ be a
twisted sum of Banach spaces, precisely \[ 0\rightarrow
H\stackrel{j}{\rightarrow} H\cplus \ E
\stackrel{q}{\rightarrow} E\rightarrow 0. \] Then \[
0\rightarrow E^* \stackrel{j^*}{\rightarrow} (H\cplus \ E)^*
\stackrel{q^*}{\rightarrow} H^*\rightarrow 0 \] is also a
twisted sum, say $E^*\cplus \ H^*$; we have \[ spl(H\cplus \
E) \ge spl(E^*\cplus \ H^*).  \] Moreover, if $E\cplus \ H$
is reflexive, equivalently if $E$ and $H$ are reflexive,
then \[ spl(H\cplus \ E) = spl(E^*\cplus \ H^*).  \]
\end{lemma}

Next we recall the profound result due to Enflo,
Lindenstrauss and Pisier \cite{ELPi} (cf.\ \cite{K},
\cite{KPc} for further examples).

\vspace{.3cm}

\noindent{\bf (ELP)} {\em There exists a twisted sum
$\ell^2\cplus \ \ell^2$ which does not split.}

\vspace{.3cm}

A simple and known consequence of (ELP) are the next two
Corollaries (cf.\ \cite[Theorem IV.6.1]{Do2}.

\begin{corollary} \label{cor-4.1} If a Banach space $E$
contains $\ell^2_n$ uniformly isomorphic and uniformly
complemented $(n=1,2,\dots)$ then there exists a twisted sum
$\ell^2\cplus \ E$ which does not split.  \end{corollary}

\bproof By (ELP) and Theorem \ref{thm-3.1} for $n=1,2,\dots$
there exists a twisted sum $Y_n=\ell^2_n\cplus \ \ell^2$
such that $\sup_n spl(Y_n)=+\infty$.  Thus, by Lemma
\ref{lemma-4.1} we have $\sup_n spl(Y_n^*)=+\infty$ because
$Y_n^*=\ell^2\cplus \ \ell^n_2$ is reflexive
($n=1,2,\dots$).  By our assumption on $E$ there is a $C\in
[1,\infty)$ such that for $n=1,2,\dots$ there are linear
operators $\varphi_n:  E\to \ell_n^2$ and $\psi_n:
\ell^2_n\to E$ satisfying (\ref{eq-4.1}) with $F=\ell_n^2$,
$\varphi=\varphi_n$ and $\psi=\psi_n$.  Thus, by Lemma
\ref{lemma-4.1} there exists a sequence $(Z_n)$ of twisted
sums $\ell^2\cplus \ E$ such that $\sup_n spl(Z_n)=+\infty$.
Thus, by Proposition \ref{prop-3.3} there exists a twisted
sum $\ell^2\cplus \ E$ which does not split.  \eproof

\begin{corollary} \label{cor-4.2} If a Banach space $E$ is
$K$-convex and infinite dimensional, then there exists a
twisted sum $\ell^2\cplus \ E$ which does not split.
\end{corollary}

\bproof By \cite[Theorem 19.3]{DJTg} every infinite
dimensional $K$-convex Banach space satisfies the assumption
of Corollary \ref{cor-4.1}.  (For the definition of
$K$-convex Banach spaces see \cite[Chapt.\ 13]{DJTg}.
\eproof

A less obvious although also ``formal'' consequence of (ELP)
is the following

\begin{theorem} \label{thm-4.1} If $E$ is an infinite
dimensional Banach space such that $E^*$ is of cotype 2,
then there exists a twisted sum $\ell^2\cplus \ E$ which
does not split.  \end{theorem}

\bproof Similarly, as in the proof of Corollary
\ref{cor-4.2}, pick for $n=1,2,\dots$ twisted sums
$Y_n=\ell_n^2\cplus \ \ell^2$ so that $\sup_n
spl(Y_n)=\sup_n spl(Y^*_n)=+\infty$.  By the dual version of
Dvoretzky's Theorem (\cite[Theorem 7.1]{To}) there exists a
bounded linear surjection $u_n:  E\to\ell_n^2$ with
$\|u_n\|[\rho(u_n)]^{-1}\le 2$.  Thus, by Lemma
\ref{lemma-1.3} there exist a twisted sum $W_n=\ell^2\cplus
\ E$ and a surjection $U_n:  W_n\to Y_n^*$ such that
$q_nU_n=u_nQ_n$ where $q_n:  Y_n^*\to\ell_n^2$ and $Q_n:
W_n\to E$ are quotient maps with the kernels isometrically
isomorphic to $\ell^2$ and $\|U_n\|\rho(U_n)\le 2$.

By Proposition \ref{prop-3.3} it is enough to show that
$\sup_n spl \, W_n=+\infty$.  Assume on the contrary that
$\sup_n spl \, W_n=C<+\infty$.  Then, by Lemma
\ref{lemma-4.2} $\sup_n spl(W_n^*)\le C$.  Thus $W_n^*$
would be $C$-isomorphic to the cartesian product $E^*\oplus
\ell^2$ hence the cotype 2 constant of $W_n^*$ would be
bounded by a constant $C_2$ independent of $n$.  Since $Y_n$
is reflexive and $U_n$ is a surjection, $U_n^*$ is an
isomorphic embedding of $Y_n$ into $W_n^*$, moreover the
Banach Mazur distance $d(Y_n,U^*_n(Y_n))\le 2$.  Thus the
cotype 2 constants of the $Y_n$'s would be uniformly bounded
by $2C_2$.  Since the $Y_n$'s are twisted sums of Hilbert
spaces, they have for every $\varepsilon>0$ uniformly
bounded type $(2-\varepsilon)$ constants (cf.\ \cite{ELPi}).

Thus, by a result of Pisier \cite{Pi}, the $K$-convexity
constants of the $Y_n$'s are uniformly bounded, say by
$C_3$.  Thus, by the Maurey-Pisier duality theorem (cf.\
\cite{MaPi}, \cite[Proposition 12.8]{To}) the $Y_n^*$'s
would have uniformly bounded type 2 constants, say by
$C_4=C_4(C_1,C_2)$.  Now let $q_n:  Y^*_n\to \ell_n^2$ be
the quotient map with $X_n=\ker q_n$ isometrically
isomorphic to $\ell^2$.  Then the map $\id_{X_n}$ could be
regarded as an operator from a subspace $X_n$ of a space
$Y_n^*$ of type 2 into a space $\ell^2$ of cotype 2. Thus,
by Maurey's extension theorem (cf.\ \cite {Ma2};
\cite[Theorem 13.13]{To}), $\id_{X_n}$ would extend to a
projection $p_n:  Y_n^*\to X_n$.  Moreover, the norms
$\|p_n\|$ would be uniformly bounded by a constant depending
on $C_4$ only.  Hence $\sup_n spl(Y_n^*)<+\infty$, a
contradiction.  \eproof

The analysis of the proofs of Corollary \ref{cor-4.2} and
Theorem \ref{thm-4.1} shows that the ``local version'' of
these assertions also holds.  Precisely one has

\begin{corollary} \label{cor-4.3} For each $C\in [1,\infty)$
there exists a sequence $(a_n^c)$ of positive numbers with
$\lim_n a_n^c=+\infty$ such that if $E$ is an
$n$-dimensional Banach space such that either the
$K$-convexity constant of $E$ or the cotype 2 constant of
$E^*$ does not exceed $C$, then there exists a twisted sum
$\ell^2\cplus \ E$ with $spl(\ell^2\cplus \ E)\ge a_n^c$
($n=1,2,\dots$).  \end{corollary}

\begin{corollary} \label{cor-4.4} If a Banach space $E$
either contains $\ell_n^\infty$ uniformly isomorphic or for
some $p\in (1,\infty)$ contains $\ell_n^p$ uniformly
isomorphic and uniformly complemented ($n=1,2,\dots$), then
there exists a twisted sum $\ell^2\cplus \ E$ which does not
split.  \end{corollary}

\bproof Recall the following well known facts.  For each
fixed $p\in (1,\infty)$ the spaces $\ell_n^p$ have uniformly
bounded $K$-convexity constants and the spaces
$(\ell_n^\infty)^*=\ell_n^1$ have uniformly bounded cotype 2
constants ($n=1,2,\dots$); the space $\ell_n^\infty$ is norm
one complemented in every larger Banach space in which it is
isometrically embedded ($n=1,2,\dots$).  Now, combine
Corollary \ref{cor-4.3} with Lemma \ref{lemma-4.2} and
Proposition \ref{prop-3.3}.  \eproof

\begin{corollary} \label{cor-4.5} For every infinite set $A$
there exists a twisted sum $\ell^2\cplus \ c_0(A)$ which
does not split.  \end{corollary}

\vspace{.3cm}

\noindent{\bf Remarks:} \begin{enumerate} \item Corollaries
\ref{cor-4.4} and \ref{cor-4.5} can be deduced from some
results of the forthcoming paper \cite{BuK} where a
different argument is used.  \item The results of this
section indicates that there are ``few'' Banach spaces $E$
for which every twisted sum $\ell^2\cplus \ E$ splits.
Since the $\{0\}$ space is a GT space, Theorem \ref{thm-3.1}
yields that if $E$ is an ${\cal L}^1$ space, then every
twisted sum $\ell^2\cplus \ E$ splits.  A more sophisticated
example is the following.  Let $X$ be a subspace of $\ell^1$
which is isomorphic to $\ell^1$ but uncomplemented in
$\ell^1$ (The existence of $X$ follows from a result of
Bourgain \cite{B}).  Put $E=\ell^1/X$.  Since $X$ being
isomorphic to $\ell^1$ is a GT-space, Theorem \ref{thm-3.1}
yields that every twisted sum $\ell^2\cplus \ E$ splits.  On
the other hand, $E$ is not an ${\cal L}_1$-space because
otherwise it would follow from Proposition \ref{prop-2.1}
that $\id_E$ lifts to $\ell^1$, i.e.\ there exists
$\widetilde{T}:  E\to\ell^1$ such that $Q\widetilde{T}=
\id_E$ where $Q:  \ell^1\to\ell^1/X$ is the quotient map.
This would yield that $X$ is complemented in $\ell^1$, a
contradiction. It seems to be an interesting problem to characterize
all Banach spaces
$E$ such that every twisted sum $\ell^2\cplus E$ splits.
\end{enumerate}
\eremark

\section{An application to Sidon sets} \label{sec-5}

\setcounter{equation}{0}

Let $G$ be a compact abelian group, $\Gamma$ its dual.
$L^p(G)$ denotes the $L^p$ space with respect to the
normalized Haar measure of $G$ denoted either by $dx$ or by
$\lambda$ ($1\le p\le\infty$).  $M(G)$ stands for the space
of all complex Borel measures on $G$ with finite variation.

Given a set $A\subset\Gamma$, put $\tilde{A}=\Gamma\backslash
A$, \[ L^p_{\tilde{A}}=\left\{ f\in L^p(G):
\hat{f}(\gamma)=0\quad\mbox{for $\gamma\in A$}\right\}
\quad(1\le p\le\infty), \] where \[ \hat{f}(\gamma) = \int_G
f(x)\overline{\gamma(x)}dx\quad (\gamma\in\Gamma,\ f\in
L^1(G)).  \]

Note that for $1\le p<\infty$, \[ L^p_{\tilde{A}}= \ \mbox{the
closed linear subspace of $L^p$ generated by $A$}.  \]

Recall that an $S\subset\Gamma$ is a {\em Sidon set} if $S$
regarded as a subspace of $L^\infty(G)$ is equivalent to the
unit vector basis of $\ell^1(S)$.  A classical
characterization of a Sidon set says

\vspace{.3cm}

\noindent{\bf (*)} {\em $S$ is a Sidon set iff the map $Q:
L^1(G)\to c_0(S)$ defined by
$Qf=(\hat{f}(\gamma))_{\gamma\in S}$ is a surjection.
Clearly $\ker Q=L^1_{\tilde{S}}(G)$.}

\vspace{.3cm}

Thus,

\begin{corollary} \label{cor-5.1} Let $S$ be an infinite
Sidon set.  Then \begin{itemize} \item[(i)] The canonical
image of $L^1_{\tilde{S}}(G)$ is uncomplemented in
$(L^1_{\tilde{S}}(G))^{**}$.  \item[(ii)] $L^1_{\tilde{S}}(G)$ not
isomorphic to a complemented subspace of a Banach lattice.
\item[(iii)] There exists a bounded linear operator from
$L^1_{\tilde{S}}(G)$ into a Hilbert space which is not
2-absolutely summing.  \item[(iv)] $L^1_{\tilde{S}}(G)$ is not
an ${\cal L}_1$-space.  \end{itemize} \end{corollary}

\bproof (i).  Combine (*) with Proposition 2.2.

(ii).  First note that $L^1_{\tilde{S}}(G)$ is a weakly
sequentially complete Banach space being a closed subspace
of a weakly sequentially complete Banach space $L^1(G)$.
Next, observe that if a Banach space $X$ is complemented in
its second dual, then every complemented subspace of $X$ has
the same property.  Furthermore, a weakly sequentially
complete complemented subspace of a Banach lattice is
isomorphic to a complemented subspace of a weakly
sequentially complete Banach lattice (\cite[II, Proposition
1.c.6]{LT}; \cite{FJT}).  Finally, a weakly sequentially
complete Banach lattice is complemented in its second dual
(\cite[II, Theorem 1.c.4]{LT}).  Combining these facts with
(i), we get (ii). Note that this argument also applies to the kernel of
any surjection onto $c_0$ in place of $L^1_{\tilde{S}}(G).$

(iii).  Combine (*) with Theorem \ref{thm-3.1} and Corollary
\ref{cor-4.5}.

(iv).  Combine (iii) with (G) in Section \ref{sec-3}.
\eproof

It is interesting to compare Corollary \ref{cor-5.1} (iii)
with the following known fact (cf.\ \cite[Theorem 2.1 and
Proposition 3.2]{KwP}).

\begin{corollary} \label{cor-5.2} Let $S\subset\Gamma$ be a
Sidon set.  Let $u:  L^1_{\tilde{S}}(G)\to H$ be a bounded
linear operator into a Hilbert space $H$.  Let $I_{\tilde{S}}:
L^2_{\tilde{S}}(G)\to L^1_{\tilde{S}}(G)$ be the natural
embedding.  Then a) $u I_{\tilde{S}}:  L^1_{\tilde{S}}(G)\to H$
belongs to the Hilbert-Schmidt class; b) if $v:
L^1_{\tilde{S}}(G)\to L^2_{\tilde{S}}(G)$ is a translation
invariant bounded linear operator, then $v$ is
$p$-absolutely summing for all $p\in [0,2]$.
\end{corollary}

\bproof Since $S$ is a Sidon set, there exists a $\mu\in
M(G)$ such that \begin{equation} \label{eq-5.0}
|\hat{\mu}(\gamma)|\ge 1 \quad\mbox{for $\gamma\in \tilde{S}$
and $\hat{\mu}(\gamma)=0$ for $\gamma\in S$.} \end{equation}
The existence of a $\mu\in M(G)$ satisfying (\ref{eq-5.0})
is a consequence of Drury's Lemma \cite{Dr} (cf.\ also
\cite{Gl} and \cite[Chapt.\ 2]{GMcG}).

For a) consider the factorization $I_{\tilde{S}}= M_\mu
Ii\hat{\mu}^{-1}_{\tilde{S}}$, i.e.  \[ I_{\tilde{S}}:
L^2_{\tilde{S}}(G)
\stackrel{\hat{\mu}_{\tilde{S}}^{-1}}{\longrightarrow}
L^2_{\tilde{S}}(G) \stackrel{i}{\longrightarrow} L^2(G)
\stackrel{I}{\longrightarrow} L^1(G)
\stackrel{M_\mu}{\longrightarrow} L_{\tilde{S}}^1(G) \] where
$\hat{\mu}_{\tilde{S}}^{-1}$ is the diagonal operator in the
character basis $(\gamma)_{\gamma\in \tilde{S}}$ taking
$\gamma$ into $\hat{\mu}(\gamma)^{-1}\cdot\gamma$ for
$\gamma\in \tilde{S}$; $i$ the natural isometrically isomorphic
inclusion; $I$ the natural embedding; $M_\mu$ the operator
of convolution with $\mu$.  Thus $uI_{\tilde{S}}$ factors
through $L^1(G)$ and the desired conclusion follows from (G)
of Section \ref{sec-3}.

For b), note that the assumption on $v$ implies that
$v(\gamma)=\alpha_\gamma\gamma$ and for $\gamma\in \tilde{S}$
for some scalar function $(\alpha_\gamma)_{\gamma\in
\tilde{S}}$.  By a), $vI_{\tilde{S}}$ is in the Hilbert Schmidt
class, hence $\sum_{\gamma\in \tilde{S}}
|\alpha_\gamma|^2=\sum_{\gamma\in \tilde{S}}
\|v(\gamma)\|^2_2<+\infty$.  Thus $\tilde{v}:  L^1(G)\to
L^2_{\tilde{S}}(G)$ defined by
$\tilde{v}(\gamma)=\alpha_\gamma\gamma$ for $\gamma\in
\tilde{S}$ and $\tilde{v}(\gamma)=0$ for $\gamma\in S$ is a
bounded linear operator and we have the factorization
$v=\tilde{v}i$ where $i:  L^1_{\tilde{S}}(G)\to L^1(G)$ denotes
the natural isometrically isomorphic inclusion.  The desired
conclusion again follows from (G) of Section \ref{sec-3}.
\eproof

Our last result shows that if $S$ is a Sidon set then the
space $L^1_{\tilde{S}}(G)$ shares an important property of
$L^1$ spaces.

Recall that a linear operator $T:  X\to Y$ is called
Dunford-Pettis or completely continuous if it takes weakly
compact sets in $X$ into norm compact sets.  A Banach space
$X$ is said to have the Dunford-Pettis property provided
every weakly compact operator from $X$ is Dunford-Pettis.

\begin{theorem} \label{thm-5.1} Let $S$ be a Sidon set.
Then $L^1_{\tilde{S}}(G)$ has the Dunford-Pettis property.
\end{theorem}

\bremark For the special case where $S$ consists of the
Rademachers, this is stated without proof by Bourgain
\cite{B}.  \eremark

The proof is based upon some properties of the Lions-Peetre
$K$-functional interpolating between the $\ell^1$ and
$\ell^2$ norms on $\Gamma$.

Let $\Psi\in L^2(G)$.  Let $(\gamma_m)$ be an enumeration of
$supp \, \widehat{\Psi}=\{\gamma\in\Gamma:
\widehat{\Psi}(\gamma)\ne 0\}$ into a sequence such that the
sequence $(|\widehat{\Psi}(\gamma_m)|)$ is non-increasing.
Put \begin{equation} \label{eq-5.1}
K_{1,2}(\widehat{\Psi},t) = \sum_{1\le m\le t^2}
|\widehat{\Psi}(\gamma_m)|+t
\left(\sum_{m>t^2}|\widehat{\Psi}(\gamma_m)|^2\right)^\frac12
\quad\mbox{for $t>0$}.  \end{equation} If the set $supp \,
\widehat{\Psi}$ is finite, say it has $m_0$ elements, then
the right hand side of (\ref{eq-5.1}) is understood to be
equal $\sum_{m=1}^{m_0} |\widehat{\Psi}(\gamma_m)|$ whenever
$t^2\ge m_0$.

The qualitative version of the following deep result of
Asmar and Montgomery-Smith (\cite[Theorem 3.9]{AM}) plays
the crucial role in our proof

\vspace{.3cm}

\noindent{\bf (AM)} {\em Let $S\subset\Gamma$ be a Sidon
set.  Then there is a constant $c>0$ depending only on $S$ so
that if $\Psi\in L^2_S(G)$, then $\lambda\{|\Psi|\ge c^{-1}
K_{1,2}(\widehat{\Psi},t)\}\ge c^{-1}e^{-ct^2}$ for all
$t>0$.}

\vspace{.3cm}

We apply the qualitative analogue of (AM) via a ``correction
lemma'' which is stated next.

\begin{lemma} \label{lemma-5.1} Let $S\subset \Gamma$.
Assume that there exist $c>0$ and a strictly decreasing
function $a:  (0,+\infty)\to (0,+\infty)$ with
$\lim_{t\to\infty} a(t)=0$ such that \[ \lambda\{|\Psi|\ge
c^{-1} K_{1,2}(\widehat{\Psi},t)\}\ge a(t)\qquad (t>0,
\Psi\in L_S^2(G)).  \] Let $\varepsilon_0=\lim_{t\to
0}\sqrt{a(t)}$.  Then there exists a function $\delta:
(0,\varepsilon_0)\to (0,+\infty)$ with $\lim_{\eta\to
0}\delta(\eta)=0$ such that if $\phi\in L^\infty (G)$ and
$\varepsilon\in (0,\varepsilon_0)$ satisfy \begin{equation}
\label{eq-5.2}
\sum_{\gamma\in\tilde{S}}|\hat{\phi}(\gamma)|^2\le
\varepsilon^2\|\phi\|^2_\infty, \end{equation} then there
exists $\varphi\in L^\infty(G)$ such that \begin{equation}
\label{eq-5.3} \begin{array}{lrcl} (i) &
\hat{\varphi}(\gamma) &=& \hat{\phi}(\gamma) \quad\mbox{for
$\gamma\in\widetilde{S}$};\\ (ii) & \|\varphi\|_\infty &\le
& (2c+1)\|\phi\|_\infty;\\ (iii) & \|\varphi\|_2 &\le &
\|\phi\|_\infty \delta(\varepsilon).  \end{array}
\end{equation} \end{lemma}

\bproof If $\phi\in L_S^\infty(G)$ put $\varphi=0$.  Let
$\phi\notin L_S^\infty(G)$.  Put $\Psi=\sum_{\gamma\in
S}\hat{\phi}(\gamma)\gamma$.  Clearly $\Psi\in L_S^2(G)$,
$\phi\ne \Psi$ and (by (\ref{eq-5.2})) \begin{equation}
\label{eq-5.4} \|\Psi-\phi\|^2_2 =
\sum_{\gamma\in\widetilde{S}} |\hat{\phi}(\gamma)|^2\le
\varepsilon^2 \|\phi\|^2_\infty.  \end{equation}

If $\varepsilon\in (0,\varepsilon_0)$ then
$\varepsilon^2=a(t)$ for some $t\in (0,+\infty)$.  Hence
\begin{equation} \label{eq-5.5} \lambda\{|\Psi|\ge c^{-1}
K_{1,2}(\hat{\Psi},t)\}\ge \varepsilon^2.  \end{equation} On
the other hand, combining (\ref{eq-5.4}) with the inclusion
\[ \{|\Psi|\ge 2\|\phi\|_\infty\}\subset \{|\Psi-\phi|\ge
\|\phi\|_\infty\} \] and with a weak type estimate \[
\lambda\{|\Psi-\phi|\ge \|\phi\|_\infty\}\le
\|\Psi-\phi\|^2_2\cdot \|\phi\|^{-2}_\infty \] we obtain
\begin{equation} \label{eq-5.6} \lambda\{|\Psi|\ge
2\|\phi\|_\infty\}\le \varepsilon^2.  \end{equation}
Moreover, the equality in (\ref{eq-5.6}) implies
$|\Psi|=2\|\phi\|_\infty 1_{\{|\Psi|\ge 2\|\phi\|_\infty\}}$
because the equality in the weak type estimate implies
$|\Psi-\phi|=1_{\{|\Psi-\phi|\ge\|\phi\|_\infty\}}$.

The inequalities (\ref{eq-5.5}), (\ref{eq-5.6}) and the
``moreover'' remark imply \begin{equation} \label{eq-5.7}
K_{1,2}(\widehat{\Psi},t)\le 2c\|\phi\|_\infty.
\end{equation}

We put $\varphi=\phi-\sum_{1\le m\le
t^2}\widehat{\Psi}(\gamma_m)\gamma_m$.  Then (\ref{eq-5.3})
(i) is obvious.  Since \[ \|\varphi\|_\infty\le
\|\phi\|_\infty +\sum_{1\le m\le t^2}
|\widehat{\Psi}(\gamma_m)|\le \|\phi\|_\infty
+K_{1,2}(\widehat{\Psi},t), \] the inequality (\ref{eq-5.7})
implies (\ref{eq-5.3}) (ii).

Let $b:  (0,\varepsilon_0)\to (0,+\infty)$ be the inverse
function of the function $\sqrt{a}$.  In particular
$b(\varepsilon)=t$.  Define $\delta:  (0,\varepsilon_0)\to
(0,+\infty)$ by \[ \delta(\eta)^2=\eta^2
+4c^2[b(\eta)]^{-2}\quad\mbox{for $\eta\in
(0,\varepsilon_0)$.} \] Clearly $\lim_{\eta\to
0}\delta(\eta)=0$ because $\lim_{\eta\to 0}
b(\eta)=+\infty$.  We have \begin{eqnarray*} \|\varphi\|^2_2
&=& \|\phi-\Psi\|^2_2 + \sum_{m>t^2}
|\hat{\phi}(\gamma_m)|^2\\ &\le &
\varepsilon^2\|\phi\|^2_\infty
+t^{-2}[K_{1,2}(\widehat{\Psi},t)]^2\quad\mbox{(by
(\ref{eq-5.1}) and (\ref{eq-5.4}))}\\ &\le &
(\varepsilon^2+4c^2[b(\varepsilon)]^{-2})\|\phi\|^2_\infty
\quad\mbox{(by (\ref{eq-5.7}))}\\ &=&
(\delta(\varepsilon)\|\phi\|_\infty)^2 \end{eqnarray*} which
verifies (\ref{eq-5.3}) (iii).  \eproof

Let $B_{\tilde{S},2}$ denote the closed unit ball of
$L^2_{\tilde{S}}(G)$.  Since $B_{\tilde{S},2}$ is a weakly
compact subset of $L^1_{\tilde{S}}(G)$, every Dunford-Pettis
operator from $L^1_{\tilde{S}}(G)$ maps $B_{\tilde{S},2}$
into a norm compact set.  Conversely one has

\begin{lemma} \label{lemma-5.2} Let $S\subset\Gamma$ satisfy
the assertion of Lemma \ref{lemma-5.1}.  Let $T:
L^1_{\tilde{S}}(G)\to Y$ ($Y$ an arbitrary normed space) be
a bounded linear operator such that $T(B_{\tilde{S},2})$ is
norm compact.  Then $T$ is Dunford-Pettis.  \end{lemma}

\bproof Let $K\subset L^1_{\tilde{S}}(G)$ be a weakly
compact set.  Then \[ \sup\{\|h\|_1:  h\in K\}=A<+\infty, \]
and $K$ is uniformly absolutely integrable, in particular,
there is a function $\omega:  (0,+\infty)\to (0,+\infty)$
such that \begin{equation} \label{eq-5.8} \begin{array}{rcl}
\lim_{\alpha\to\infty} \omega(\alpha) &=&0,\\
\|h-h^\alpha\|_1 &\le & \omega(\alpha) \quad\mbox{for
$\alpha>0$ and $h\in K$}, \end{array} \end{equation} where
$h^\alpha = h\cdot 1_{\{|h|\le\alpha\}}$.

We shall show that $T(K)$ is a norm totally bounded set,
hence it is norm compact.  To this end it suffices to verify
\begin{equation} \label{eq-5.9} \begin{array}{l} \mbox{for
every $\sigma>0$ there exists $\varepsilon>0$ such that}\\
T(K) \subset \sigma B_Y \cup T\left(\frac{2A}{\varepsilon}
B_{\tilde{S},2}\right) \end{array} \end{equation} ($B_Y$
denotes the unit ball of $Y$).

Fix $f\in K$.  Let $\varepsilon\in (0,\varepsilon_0)$ be
given ($\varepsilon_0$ is that of Lemma \ref{lemma-5.1}).
We shall choose $\varepsilon$ for $\sigma$ later on.  Put \[
\rho = \inf \{\|T(f)-T\left(\frac{2A}{\varepsilon}
g\right)\|:  g\in B_{\tilde{S},2}\}.  \] By the Separation
Theorem (\cite[Theorem V.2.12]{DnS}) there exists a linear
functional $x^*\in Y^*$ such that \[ \rho = \sup\{\Re \,
x^*(y):  \|y\|\le\rho\} = \inf\{\Re \,
x^*(T(f)-T\left(\frac{2A}{\varepsilon} g\right)):  g\in
B_{\tilde{S},2}\}.  \] Clearly $\|x^*\|=1 $. Let
$\phi^*=T^*(x^*)\in (L^1_{\tilde{S}}(G))^*$ and let $\phi\in
L^\infty(G)=(L^1(G))^*$ be a norm preserving extension of
$\phi^*$.  Then \begin{eqnarray*} \int_G \phi h \, dx &=&
x^*(T(h)) \quad\mbox{for $h\in L^1_{\tilde{S}}(G)$},\\
\|\phi\|_\infty &\le & \|T\|.  \end{eqnarray*} Since $0\in
B_{\tilde{S},2}$, we have \begin{equation} \label{eq-5.10}
\rho \le |\int_G\phi f\, dx| \le A\|\phi\|_\infty.
\end{equation} Thus, if $g\in B_{\tilde{S},2}$ then \[ \Re
\int \phi \frac{2A}{\varepsilon} g\, dx \ge -\rho +\Re \int
\phi f\, dx \ge -2A\|\phi\|_\infty.  \] Consequently, taking
into account that $B_{\tilde{S},2}$ is a circled set \[
|\int\phi g\, dx|\le \varepsilon
\|\phi\|_\infty\quad\mbox{for $g\in B_{\tilde{S},2}$.} \]
The latter inequality implies (\ref{eq-5.2}).

Now, given $\sigma>0$ invoking (\ref{eq-5.8}), one picks
$\alpha>0$ so that \[ \|f-f^a\|_1\le\sigma
(2(2c+1)(\|T\|+1))^{-1}\quad\mbox{for $f\in K$.} \] Next we
choose $\varepsilon\in (0,\varepsilon_0)$ so that
$\delta(\varepsilon)\le
\sigma(2\sqrt{\alpha}(\|T\|+1))^{-1}$ where $\delta(\cdot)$
is the function of Lemma \ref{lemma-5.1}.  Finally, for
$\phi$ which has been constructed for fixed $f\in K$ and
$\varepsilon$ just chosen ($\phi$ satisfies (\ref{eq-5.2})
and $\|\phi\|_\infty\le\|T\|$) we apply Lemma
\ref{lemma-5.1} to construct $\varphi$ satisfying
(\ref{eq-5.3}).  Then \begin{eqnarray*} |\int\phi f \, dx|
&=& |\int\varphi f\, dx| \qquad\mbox{(by (i) because $f\in
K\subset L^1_{\tilde{S}}(G)$)}\\ &\le & |\int_G
\varphi(f-f^\alpha)dx|+|\int_G\varphi f^\alpha \, dx|\\ &\le
& \|\varphi\|_\infty
\|f-f^\alpha\|_1+\|\varphi\|_2\sqrt{\alpha}\\ &\le & \sigma
\qquad\mbox{(by (ii) and (iii))}.  \end{eqnarray*} Thus, by
(\ref{eq-5.10}), $\rho\le\sigma$.  Hence $T(f)\in\sigma
B_Y\cup T(\frac{2A}{\varepsilon}B_{\tilde{S},2})$ which
yields (\ref{eq-5.8}).  \eproof

\begin{lemma} \label{lemma-5.3} Assume that for some
$S\subset\Gamma$ there exists a $\mu\in M(G)$ satisfying
(\ref{eq-5.0}).  Let $T$ be a weakly compact operator from
$L^1_{\tilde{S}}(G)$.  Then $T(B_{\tilde{S},2})$ is a norm
compact set.  \end{lemma}

\bremark Note that the condition ``$T(B_{\tilde{S},2})$ norm
compact'' is equivalent to ``the operator
$TI_{|L^2_{\tilde{S}}(G)}$ is compact'' where $I:  L^2(G)\to
L^1(G)$ is the natural embedding.  \eremark

\bproof Let $M_\mu:  L^1(G)\to L^1_{\tilde{S}}(G)$ be the
operator of convolution with $\mu$ (by (\ref{eq-5.0}),
$M_\mu(L^1(G))\subset L^1_{\tilde{S}}(G)$).  Then the set
$TM_\mu(B_{\tilde{S},2})$ is norm compact because $L^1(G)$
has the Dunford-Pettis property (\cite[Theorem
VI.8.12]{DnS}).  Then $T(B_{\tilde{S},2})$ is norm compact
because if $M_\mu$ is regarded as an operator from $L^2(G)$
into $L^2_{\tilde{S}}(G)$ then the restriction
$M_{\mu|L^2_{\tilde{S}}(G)}$ is, by (\ref{eq-5.0}),
invertible.  \eproof

\noindent{\bf Proof of Theorem \ref{thm-5.1}.} If
$S\subset\Gamma$ is a Sidon set, then, by \cite{AM}, $S$
satisfies the assumption of Lemma \ref{lemma-5.1}, while by
the result of Drury \cite{Dr} mentioned above, $S$ satisfies
the assumption of Lemma \ref{lemma-5.3}.  Now Theorem
\ref{thm-5.1} follows directly from Lemmas
\ref{lemma-5.1}--\ref{lemma-5.3}.  \eproof

 \bremark We do not know if a subspace $X$ of $L_1$ has the
Dunford-Pettis property whenever $L_1/X$ is isomorphic to $c_0.$\eremark

\noindent Department of Mathematics, University of Missouri
Columbia, Columbia, MO 65211, USA\\ e-mail
MATHNJK@MIZZOU1.MISSOURI.EDU\\
 \ \\ Institute of Mathematics, Polish Academy of Sciences,
Sniadeckich 8, Ip.\ 00950 Warszawa, POLAND\\ and \\
Department of Mathematics, University of Missouri Columbia,
Columbia, MO 65211, USA\\ e-mail OLEK@IMPAN.IMPAN.GOV.PL

\end{document}